\documentclass[12pt]{article} 
\usepackage{amsfonts}
\usepackage{amsmath}
\usepackage{graphicx}
\usepackage{epsfig}
\usepackage{amssymb}
\usepackage{amstext}

\newtheorem{theorem}{Theorem}[section]
\newtheorem{lemma}{Lemma}[section]
\newtheorem{remark}{Remark}[section]
\newtheorem{corollary}{Corollary}[section]

\newcommand\bcoper{\vartriangleright{\kern-.55em\hbox{$\blacktriangleleft \,
$}}}

\def\eendproof{\ \ \hfill{$\Box$}}

\begin{document}



\title{Billiards and the Five Distance Theorem II}


\author{Jan~Florek}
\date{}

\maketitle


\section{Introduction}
This paper is a complement to the author earlier paper \cite{F}.
We consider a billiard table rectangle with perimeter of length $1$. If a billiard ball is sent out from position $F(1)$ at the angle of $\pi/4$, then the ball will rebound against the sides of the rectangle consecutively in points $F(2)$, $F(3)$, \dots.
Let $n\geq5$ and $\Phi= \{ F(j): 1\leq j\leq n \}$ be the set of different points.
An open connected subset of the perimeter of the billiard rectangle with different endpoints from the set $\Phi$ is called \textit{segment}. A segment with endpoints $F(k)$, $F(l)$, $1\le k,l\le n$, is called \textit{even} (or \textit{odd}), and has \textit{weight} $|k-l|$ (or $k+l$) if $k$, $l$ are of the same (or different) parity.  \textit{Length} of a segment is a distance along the perimeter between its endpoints. A segment with the length less than $\frac{1}{2}$ is called \textit{short}. A segment is called \textit{elementary} if there are no points of the set $\Phi$ between its endpoints. Suppose $\emptyset \neq V\subseteq\{F(1),F(n)\}$. A segment $I$ is \textit{associated} with $V$ if $I$ is an elementary segment incident with an element of $V$ or $I\cap\Phi$ is nonempty set contained in $V$. Since $n\geq5$, any segment associated with $\{F(1),F(n)\}$ is determined by its endpoints uniquely. By $(F(k),F(l))$ (respectively $(F(k),F(l))_e$) we denote the short segment (respectively the elementary segment) incident with $F(k)$ and $F(l)$. Note that $(F(k),F(l))$ and  $(F(l),F(k))$ are the same segments.

Let $\omega_1<\omega_2$ be odd weights and $\omega_0$
be an even weight of segments associated with $\{F(1)\}$, and let $\omega_3<\omega_4$ be other odd weights of segments
associated with $\{F(1),F(n)\}$ (see Corollary~\ref{corr2.1}). Let $a_i$ be the length of the segment with the weight $\omega_i$, and suppose that $A_i$ is the set of all elementary segments with weight $\omega_i$, $i = 0, \ldots, 4$.  The author \cite{F} have proved that the weights of elementary segments have at most five different values $\omega_0, \ldots, \omega_4$. Hence,
$$
|A_0|+ \cdots + |A_4|=n.
$$
If $a_2-\epsilon a_1 = a_3-\delta a_4 =a_0$, for some $\epsilon$, $\delta\in\{-1,1\}$, then the following equations are satisfied
$$
|A_2|+\epsilon|A_1|=|A_3|+\delta|A_4|=\frac{1}{2}\omega_0.
$$
Moreover, elementary segments with equal weights have equal lengths.
In this paper we prove the following equalities (see Theorem~\ref{theo2.2}):
$$
|A_1|= \frac{\omega_1 - 1}{2}, \hbox{ and } |A_4|= n - \frac{\omega_4 - 1}{2}.
$$

Notice that we can consider a general case, when the initial angle of the ball's motion is not $\pi/4$.
By a linear transformation we can change the billiard table rectangle,
so that the general case is transformed to the $\pi/4$ case (except  the result that elementary segments with equal weights have equal lengths).
\section{The main result}
Recall the following results of the paper \cite{F}.
\begin{remark}\label{rema2.1} \cite[Remark 2.2]{F}
An odd elementary segment which is not short is of the form $(F(1), F(2))_e$ or $(F(n-1), F(n))_e$.
\end{remark}
\begin{theorem}\label{theo2.1}\cite[Theorem 2.2]{F}
If  $(F(k), F(1))_e$, $(F(1), F(n))_e$, $1<k<n$, are elementary segments incident with $F(1)$, then $k$ is even.
\end{theorem}
\begin{corollary}\label{corr2.1}\cite[Corollary 2.1]{F}
If $(F(1), F(n))_e$, $(F(n), F(l))_e$, $1<l<n$ are elementary segments incident with $F(n)$, then $l, n$ are of different parity.
\end{corollary}
\begin{lemma}\label{lemma2.1}\cite[Lemma 2.3(1)]{F}
Let $k+1<l$ and $k$, $l$ be of different parity. Then we have:
 a short segment $(F(k), F(l))$ is elementary, or is associated with $\{F(1),F(n)\}$
if and only if
the short segment $(F(k+1), F(l-1))$ is elementary.
\end{lemma}
\begin{corollary}\label{corr2.1}\cite[Corollary 2.2]{F}
There exist exactly two odd segments and one even segment associated with $\{F(1)\}$, and exactly four odd segments associated with $\{F(1), F(n)\}$. The weights of these four odd segments are different.
\end{corollary}
Now we present the main result of this paper.
\begin{theorem}\label{theo2.2}
Let $\omega_1<\omega_2$ be odd weights and $\omega_0$ be an even weight of segments associated with $\{F(1)\}$, and let $\omega_3<\omega_4$ be other odd weights of segments associated with $\{F(1),F(n)\}$. If $A_1$ $($$A_4$$)$ is the set of all elementary segments with weight $\omega_1$ $($$\omega_4$, respectively$)$,
then
$$
|A_1|= \frac{\omega_1 - 1}{2}\quad\hbox{and}\quad |A_4|= n - \frac{\omega_4 - 1}{2}.
$$
\end{theorem}
{\bf Proof.}
Suppose $(F(k),F(1))_e$, $(F(1),F(l))_e$, $k<l$, are elementary segments incident with $F(1)$, and $(F(r),F(n))_e$, $(F(n),F(t))_e$, $r<t$, are elementary segments incident with $F(n)$. By Theorem~2.1 and Corollary~2.1 $k$ is even and $t$, $n$ are of different parity. Hence, segments  $(F(k),F(1))_e$, $(F(n),F(t))_e$ are odd. If $(F(k),F(1))_e$ (or $(F(n),F(t))_e$) is short, then by Lemma~2.1,
$$|A_1| = \frac{k}{2}= \frac{\omega_1 -1}{2}$$
(or
$$|A_4| = \frac{n-t+1}{2}= \frac{2n-\omega_4 + 1}{2},$$
respectively). If $(F(k),F(1))_e$ (or $(F(n),F(t))_e$) is not short, then by Remark~2.1 and Lemma~2.1
$$k =2 \quad \hbox{and} \quad  |A_1|= 1 = \frac{\omega_1 -1}{2}$$
(or
$$t = n-1 \quad \hbox{and} \quad  |A_4|= 1 = \frac{2n-\omega_4 + 1}{2}$$
respectively).
 \eendproof

\noindent
{Institute of Mathematics,\\University of Economics \\ul. Komandorska 118/120\\ 53--345 Wroc{\l}aw, Poland\\
E-mail: jan.florek@ue.wroc.pl}
\end{document}